\documentclass{article}

\usepackage[T2A]{fontenc}
\usepackage[cp1251]{inputenc}
\usepackage[tbtags]{amsmath}
\usepackage{amsfonts,amssymb}

\usepackage{mathrsfs}

\usepackage{graphicx}

\newtheorem{theorem}{Theorem}
\newtheorem{proposition}{Proposition}

\newtheorem{corollary}{Corollary}

\newtheorem{problem}{Question}

\def\D{{\cal D}}
\def\R{{\mathbb R}}

\def\Z{{\mathbb Z}}
\def\N{{\mathbb N}}
\def\H{{\mathcal H}}

\def\V{{\mathcal{V}}}
\def\Q{{\mathcal{Q}}}
\def\g{{\mathfrak{g}}}

\begin{document}

\date{}

\author{I.A.~Taimanov
\thanks{Sobolev Institute of Mathematics, Novosibirsk; taimanov@math.nsc.ru}}

\title{Central extensions of Lie algebras, dynamical systems, and symplectic nilmanifolds
\thanks{The work is supported by RSCF (project 24-11-00281).}
}

\maketitle

\begin{abstract}
The connections between Euler's equations on central extensions of Lie algebras and Euler's equations on the original, extended algebras are described. A special infinite sequence of central extensions of nilpotent Lie algebras constructed from the Lie algebra of formal vector fields on the line is considered, and the orbits of coadjoint representations for these algebras are described. By using the compact nilmanifolds constructed from these algebras by I.K. Babenko and the author, it is shown that covering Lie groups for symplectic nilmanifolds can have any rank as solvable Lie groups.
\end{abstract}

\vskip3mm

\hfill{\it To S.V. Bolotin and D.V. Treshchev on the anniversaries}

\vskip3mm

\section{Euler's equations on central extensions of Lie algebras and magnetic geodesic flows}

Let $G$ be a simply connected Lie group and $\g$ be its algebra. In the sequel, we assume that all Lie algebras are real.

A one-dimensional central extension of a Lie algebra $\g$ is a Lie algebra $\widehat{\g}$ such that
the exact sequence
\begin{equation}
\label{central}
0 \to \R \to \widehat{\g} \to \g \to 0,
\end{equation}
and the subalgebra $\R$ generated by some element $w \in \widehat{\g}\setminus \g$
lies in the center of $\widehat{\g}$.
Obviously
$$
[u,v]_{\widehat{\g}} = [u,v]_{\g} + \alpha(u,v) w, \ \ u,v \in \g,
$$
where $\alpha$ is a skew-symmetric bilinear $2$-form on $\g$.
From the Jacobi identity for $\widehat{\g}$ it immediately follows that the form $\alpha$
is closed and hence realizes the cohomology class
$$
[\alpha] \in H^2(\g;\R).
$$
The cohomology group $H^2(\g;\R)$ parameterizes the equivalence classes of exact sequences of the form \eqref{central}.

Not every central extension of a Lie algebra corresponds to a central extension of the Lie group
$$
1 \to \R \to \widehat{G} \to G \to 1,
$$
however, in the case when the group $G$ is simply connected, this is certainly the case and the Lie algebra
$\widehat{\g}$ is a tangent algebra at the unit of $\widehat{G}$ \cite{Neeb}.

Let a Hamiltonian system with left-invariant Hamiltonian $H$ be defined on the cotangent space $T^\ast G$.
By symplectic reduction, it can be reduced to the Euler equations on the Lie coalgebra $\g^\ast$:
\begin{equation}
\label{euler}
\frac{df}{dt} = \{f, H\},
\end{equation}
where $\{ \cdot,\cdot\}$ are the Lie--Poisson brackets \cite{Arnold}.

The tangent and cotangent spaces at the unit of a Lie group are naturally identified
with the Lie algebra $\g$ and coalgebra $\g^ast$ of this group.

The coordinate functions $x_i, i=1,\dots, n =\dim \g$, on the Lie coalgebra $\g^\ast$ are elements of $\g$, namely, linear functions on $\g^\ast$. Since, by definition, the Poisson bracket satisfies the relation
$$
\{f,g\} = \sum_{i,j} \frac{\partial f}{\partial x_i}\frac{\partial g}{\partial x_j} \{x_i,x_j\},
$$
it suffices to define it for the basis functions: $\{x_i,x_j\}$. The Lie--Poisson bracket is defined by the relations
\begin{equation}
\label{lp}
\{x_i,x_j\} = \sum_k c^k_{ij} x_k,
\end{equation}
where $c^k_{ij}$ are the structure constants of the corresponding Lie algebra $\g$:
$$
[e_i,e_j] = \sum_k c^k_{ij} e_k.
$$
It is known that the restrictions of the Lie--Poisson brackets to the orbits of the coadjoint representation define symplectic structures on them, and these orbits are invariant with respect to the Hamiltonian flows \eqref{euler} \cite{Kir} (see also \cite{Arnold}).

The geodesic flow is determined by the Hamiltonian function quadratic in the velocities
$$
H = \frac{1}{2} \sum_{i,k} g^{ik} x_i x_k,
$$
which is the result of the symplectic reduction of the geodesic flow on  $T^\ast G$:
$$
H(q,p) = \frac{1}{2} \sum_{i,k} g^{ik}(q) p_i p_k, \ \ q \in G,
$$
where the momenta $p$ are related to the velocities by
$$
\dot{q}^i = g^{ik} p_k
$$
and the symmetric tensor $g^{ik}(q)$  is invariant under left translations on $G$. In this case, the Poisson brackets on $T^\ast G$ are standard:
\begin{equation}
\label{poisson}
\{q^i, p_k\} = \delta^i_k, \ \ \{p_i,p_k\} = \{q^i,q^k\} = 0, \ \ 1 \leq i,k \leq n =\dim G.
\end{equation}
The magnetic geodesic flow on $G$ is determined by the same Hamiltonian function, only the Poisson brackets \eqref{poisson} are deformed: on pairs of momenta coordinates, theire values are equal to
$$
\{p_i,p_k\} = B_{ik}(x),
$$
where $B$ is a closed $2$-form defining the magnetic field \cite{Novikov}.
The following fact is obvious.

\begin{proposition}
A left-invariant closed $2$-form $B$ on a Lie group $G$ defines a $2$-cocycle $\alpha_B$
on the Lie algebra $\g$.
\end{proposition}

Let us consider the central extension $\g_B$ of the Lie algebra $\g$ defined
by the cocycle $\alpha_B$. Let $x_{n+1}$ corresponds to the central element which does not lie in $\g$. We consider the Hamiltonian
\begin{equation}
\label{ham2}
\widehat{H} (x) = \frac{1}{2} \left(\sum_{i=1}^n x_i^2 + x_{n+1}^2\right) = H(x) +
\frac{1}{2}x_{n+1}^2,
\end{equation}
where
$x_1,\dots,x_n$ are orthonormal coordinates with respect to some scalar product on $\g$.

If the group $G$ is simply connected, then the Hamiltonian \eqref{ham2} corresponds to the
geodetic flow of the left-invariant metric on the central extension $G_B$ of the group $G$.

If the Lie algebra $\g$ is nilpotent, then such groups $G$ and $G_B$ exist.
In what follows, for brevity, we will simply speak of Euler equations for flows on
Lie algebras.

The Euler equations on $\g^\ast_B$ for the Hamiltonian \eqref{ham2}
$$
\dot{x}_i= \{x_i,\widehat{H}\}_{\g_B}, \ \ i=1,\dots,n+1,
$$
take the form
\begin{equation}
\label{eulmag}
\begin{array}{c}
\dot{x}_i = \{ x_i, H\}_{\g} + x_{n+1} \sum_{k=1}^n B_{ik} x_k, \ \ i=1,\dots,n, \\
\dot{x}_{n+1} = 0.
\end{array}
\end{equation}

\begin{proposition}
The Euler equations \eqref{eulmag} on $\g_B^\ast$, i.e. the geodesic flow of the metric \eqref{ham2},
for each fixed value of $x_{n+1}$ describes the magnetic geodesic flow on
$G_B$ corresponding to the metric $H$ and the magnetic field $x_{n+1} B$.
\end{proposition}

This proposition is obtained by directly applying symplectic reduction to magnetic geodesic flows and is used, for example, in \cite{E}.

Another observation can be made: since $\{x_{n+1},f\}_{\g_B}=0$ for any function $f$, then
$$
\dot{x}_i = \{x_i,\widehat{H}\}_{\g_B} = \{x_i,H\}_{\g_B}, \ \ i=1,\dots,n+1,
$$
and the latter system defines the Euler equations for the normal sub-Riemannian geodesic flow on the Carnot group $G_B$ with a non-integrable distribution of codimension one given by the subalgebra $\g$.

The Euler equations on Lie algebras for left-invariant normal geodesic flows of sub-Riemannian metrics were considered in \cite{MSS,T97}, and the connection of such flows with magnetic geodesic flows in a more general aspect was first established in \cite{M}.

Let's combine all these known facts as follows:

\begin{theorem}
1) Euler's equations \eqref{eulmag} of the geodesic flow on $G_B$ corresponding to the left-invariant metric \eqref{ham2} are equivalent to Euler's equations for the normal geodesic flows of the left-invariant sub-Riemannian metric given by the Hamiltonian $H$ and the non-integrable distribution generated by $\g$;

2) Euler's equations \eqref{eulmag} for each fixed value $x_{n+1} = \mathrm{const}$ define Euler's equations of the magnetic geodesic flow on $G$ given by the Hamiltonian $H$ and the magnetic field $x_{n+1} B$;

3) If Euler's equations \eqref{eulmag} are integrable, then the Euler equations for magnetic geodesic flows on $G$ are integrable for almost all $x_{n+1}$.
\end{theorem}

Note that if all structure constants of a nilpotent Lie algebra $\g$ are integer: $c^k_{ij} \in \Z$, and the cocycle
$\alpha_B$ takes integer values ??on the basis vectors
$e_1,\dots,e_n$: $[\alpha_B] \in H^2(\g;\Z)$, then the geodesic flows descend to closed nilmanifolds.

We do not dwell on the concept of integrability, but only note that $x_{n+1}$ is a first integral of any Hamiltonian system on $\g_F$ and the presence of a sufficient number of first integrals of the system \eqref{eulmag} that are functionally independent almost everywhere implies that they are functionally independent on almost all hypersurfaces $\{x_{n+1} = \mathrm{const}\}$. Item 3 of Theorem 1 shows that the question of integrability of the Euler equations on central extensions of Lie algebras deserves special attention.

\section{Note on Casimir polynomials}

A Casimir function for a Lie algebra $\g$ is, by definition, a function on $\g^\ast$ that is invariant under the coadjoint representation of $G$.

Finding Casimir functions $F$ for an algebra $\g$ with generators $e_1,\dots,e_n$
is reduced to solving systems of linear equations for the derivatives of $F$:
\begin{equation}
\label{s1}
\{x_i,F\} = \sum_{k=1}^n \frac{\partial F}{\partial x_k} \{x_i,x_k\} = \sum_{k=1}^n A_{ik} \frac{\partial F}{\partial x_k} = 0, \ \ \ i=1,\dots,n,
\end{equation}
where $\{ \cdot, \cdot\}$ is the Lie--Poisson bracket \eqref{lp} and
\begin{equation}
\label{A}
A = \left(A_{ik} = \{x_i,x_k\}\right).
\end{equation}

We will consider only Casimir functions which polynomials in $x_1,\dots,x_n$.
They form a polynomial ring $\R[u_1,\dots,u_\nu]$,
where
\begin{equation}
\label{nu}
\nu = n - \mathrm{rank}\, A,
\end{equation}
and the rank of the matrix \eqref{A} is estimated at a generic point \cite{BB}.

The upper bound for $\nu$ is obvious.
The solutions $s=(s_1,\dots,s_n)$ of the system
\begin{equation}
\label{s11}
\sum_k \{x_i,x_k\} s_k =0,
\end{equation}
to be gradients of some functions $F$, must satisfy
the compatibility conditions
$$
\frac{\partial s_k}{\partial x_i} = \frac{\partial s_i}{\partial x_k}, \ \ i,k=1,\dots,n.
$$
But there are certainly $n-\nu$ functionally independent such solutions.
The equality in \eqref{nu} is proved as follows. In the neighborhood of a typical orbit of the coadjoint representation, which has the maximum dimension, we can introduce local coordinates $y_1,\dots,y_n$, in which the orbit is defined by the equations
$$
y_{2k+1} = \dots = y_n =0,
$$
where $2k$ is the dimension of the orbit. The equations \eqref{s1} in these coordinates take the form
$$
\frac{\partial F}{\partial y_1} = \dots = \frac{\partial F}{\partial y_{2k}} = 0.
$$
Obviously, functions depending only on $y_{2k+1},\dots,y_n$ satisfy these equations.

Note that we can choose a basis of solutions to \eqref{s11} such that they are polynomi\-als in $x_1,\dots,x_n$. The Casimir functions $F$ constructed from them will also be polynomials in these variables.

\section{A tower of nilpotent Lie groups with integrable magnetic geodesic flows}

Natural examples to which Theorem 1 is applicable are given by nilpotent Lie algebras
$\Q_n$.

Recall that a Lie algebra $\g$ is called nilpotent if its lower central series
$\g_1=\g, \g_2 = [\g,\g], \dots, \g_k=[\g,\g_{k-1}], \dots$ degenerates at a finite step:
$\g_m = 0$ for some $m$.

The Campbell--Hausdorff formula
\begin{equation}
\label{ch}
u \times v = u + v + \frac{1}{2}[u,v] + \frac{1}{12} ([u,[u,v]] + [v,[v,u]]) +\dots
\end{equation}
for nilpotent algebras reduces to a finite sum and defines multiplication by $\R^n \approx \g$.
The resulting $n$-dimensional manifold with multiplication will be a Lie group $G$ with Lie algebra $\g$.

We consider an infinite series of nilpotent Lie algebras $\Q_n$ such that
$\dim \Q_n = n$ and the basis elements $e_1,\dots,e_n$
satisfy the commutation relations:
$$
[e_1,e_k] = - [e_k, e_1] = e_{k+1}, \ \ k=2,\dots,n-1,
$$
and the remaining commutators are zero.

Obviously, for any $n$, the algebra $\Q_n$ is a central extension of the algebra
$\Q_{n-1}$ obtained by adding the central element $e_n$. They correspond to nilpotent Lie groups $Q_n$, which form a tower of central extensions
$$
Q_1 = \R \subset Q_2 = \R^2 \subset Q_3 = \H_3 \subset \dots \subset Q_n \subset \dots,
$$
where
$$
\H_3 = \left\{\begin{pmatrix} 1 & x & z \\ 0 & 1 & y \\ 0 & 0 & 1 \end{pmatrix} \ : \ \ x,y,z \in \R\right\},
$$
is the three-dimensional Heisenberg group.

It is well known that typical orbits of coadjoint representations for algebras $\Q^\ast_n$ are two-dimensional (see, e.g., \cite[Corollary 1.8]{Butler}) and therefore the obvious statement holds.

\begin{proposition}
The Euler equation for any left-invariant Hamiltonian such that its value
on $T^\ast_e Q_n = \Q_n^\ast$ is not a Casimir polynomial is integrable.
\end{proposition}

\begin{corollary}
The lef--invariant magnetic geodesic flows on $Q_n$ and the normal geodesic flows of sub-Riemannian metrics on Carnot groups $Q_n$ with left-invari\-ant distributions generated by subalgebras $\Q_{n-1}$ are integrable.
\end{corollary}

The integrability of the Euler equations on the Lie algebra $\g$ implies the integrability of the corresponding Hamiltonian system on the simply connected Lie group $G$. However, this does not imply the integrability of this system dropped onto the homogeneous quotient space $G/\Gamma$, where $\Gamma$ is a group acting freely and totally discontinuously on $G$. The simplest example is provided by hyperbolic surfaces: the geodesic flow of a metric of constant negative curvature is integrable on the Lobachevsky plane $\mathbb{H}^2$, but is chaotic (Anosov) on compact quotient spaces $\mathbb{H}^2/\Gamma$.

In \cite{Butler}, the integrability of geodesic flows was proved on compact nilmanifolds $Q_n/\Gamma_n$. This proof can be extended to the case of magnetic and normal sub-Riemannian geodesic flows. It is essential that in the case of geodesic flows, when passing to nilmanifolds, it is impossible to find a complete set of analytic first integrals; one of them must necessarily be non-analytic, but it can be chosen in the class $C^\infty$.

Butler's work \cite{Butler} generalizes his earlier result \cite{Butler99,Butler00}, which proved the
$C^\infty$-integrability of the left-invariant geodesic flow on the three-dimensional nilmanifold
$\H_3/\H_3(\Z)$, where
$\H_3(\Z)$ is the subgroup of $\H_3$ formed by all elements with $x,y,z \in \Z$. Clearly,
$\H_3 = Q_3$.
The fundamental group of the manifold $\H_3/\H_3(\Z)$ is not almost commutative, which, according to the author's results, prevents complete integrability in terms of analytic first integrals \cite{T87,T88} and shows that the topological obstructions to analytic integrability obtained in these papers do not generalize to the case of $C^\infty$-integrability.

The nilmanifold $\H_3/\H_3(\Z)$ can be viewed from another point of view: as a superstructure of an automorphism of a two-dimensional torus:
$$
\H_3/\H_3(\Z) = (\R^2/\Z^2 \times [0,1]) /\{(X,0) \sim (CX,1)\}, \ \
C = \begin{pmatrix} 1& 1 \\ 0 & 1 \end{pmatrix}, \ \ X = \begin{pmatrix} x \\ y \end{pmatrix}.
$$
In \cite{BT} (see also \cite{BT2}) Butler's trick for constructing $C^\infty$-first integrals was generali\-zed to the case of the hyperbolic automorphism of the torus
$$
C = \begin{pmatrix} 2& 1 \\ 1 & 1 \end{pmatrix}.
$$
The solvmanifold constructed gave the first known example of a closed Rieman\-nian manifold
(with an analytic Riemannian metric) on which the geodesic flow is $C^\infty$-integrable but has positive topological entropy.

The theorem on topological obstructions from \cite{T87} is a multidimensional generalization of Kozlov's theorem that on a two-dimensional oriented closed surface, geodesic flows
can be analytically integrable only on the sphere and torus \cite{Kozlov} (see also
\cite{Kozlov83}). In this connection, we note that it remains open

\begin{problem}
Do there exist $C^\infty$-integrable
geodesic flows of Riemannian or Finsler metrics on closed oriented two-dimensional surfaces?
\end{problem}

\section{the Lie algebras $\V_n$, their Casimir polynomials and the orbits of the coadjoint representation}
\label{secvn}

In \cite{BaT} an infinite series of nilmanifolds with interesting properties was constructed.
Their construction was based on the infinite-dimensional Witt algebra $W(1)$ of formal vector
fields on the line with generators
$$
e_k = x^{k+1}\frac{d}{dx}, \ \ k=-1,0,1,\dots,
$$
for which the commutation relations
$$
[e_i, e_j] = (j-i) e_{i+j}, \ \ i,j \geq -1
$$
hold.

The algebra $W(1)$ has a natural filtration
$$
\dots \subset \mathcal{L}_1(1) \subset \mathcal{L}_0(1) \subset \mathcal{L}_{-1}(1) = W(1),
$$
where $\mathcal{L}_k(1), k \geq -1$, denotes the Lie subalgebra generated by the elements
$e_k, e_{k+1}, \dots$.
The quotient Lie algebras $\V_n$ defined as
$$
\V_n = \mathcal{L}_1(1)/\mathcal{L}_{n+1}(1)
$$
are nilpotent Lie algebras with generators
$e_1,\dots$,$e_n$ satisfying the commutation relations
$$
[e_i,e_j] = \begin{cases} (j-i)e_{i+j}, & i+j \leq n; \\ 0, & i+j>n. \end{cases}
$$
By means of the Campbell--Hausdorff formula, the nilpotent Lie groups $V_n$ are constructed from the Lie algebras $\V_n$, which are used to define the compact nilmanifolds $V_n$ introduced in \cite{BaT}.
We will discuss some of them in \S \ref{secnm}, and now we will move on to describing the orbits of coadjoint representations in coalgebras $\V_n^\ast$.

Casimir polynomials are invariants of the coadjoint representation.
To find Casimir polynomials, it is necessary to substitute into \eqref{s1} the expressions for
the values of the Lie--Poisson brackets on the coordinate functions:
$$
\{x_i,x_j\} = \begin{cases} (j-i)x_{i+j}, & i+j \leq n; \\ 0, & i+j>n. \end{cases}
$$
Obviously, the (coordinate) function $x_n$ is a Casimir function for the algebra $\V_n$ for any $n$.

The matrix $A$ for $\V^\ast_n$ takes the form
$$
\begin{pmatrix}
0 & x_3 & 2x_4 & \dots & (n-3) x_{n-1} & (n-2) x_n & 0 \\
-x_3 & 0 & x_5 &  \dots & (n-4) x_n & 0 & 0 \\
\dots & \dots & \dots & \dots & \dots & \dots & \dots \\
-(n-3) x_{n-1} & -(n-4) x_n & 0 & 0 & 0 & 0 & 0 \\
-(n-2) x_n & 0 & 0 & 0 & 0 &  0 & 0 \\
0 & 0 & 0 & 0 & 0 & 0 & 0
\end{pmatrix} =
$$
\begin{equation}
\label{Av}
=
\begin{pmatrix}
& & & 0\\
& \widetilde{A} & & 0 \\
& & & 0 \\
0 & \dots & 0 & 0
\end{pmatrix}
\end{equation}

\begin{proposition}
The rank of the matrix \eqref{Av} at a generic point depends on the parity of  $n$ and

1) for $n=2q+1$ and $x_n \neq 0$ the rank of the matrix \eqref{Av} is equal to $n-1$;

2) for $n=2q+2$ and $x_n \neq 0$ the rank of the matrix \eqref{Av} is equal to $n-2$.
\end{proposition}

{\sc Proof.} All entries that are to the right of the ``anti-diagonal''
$\widetilde{A}_{i,n-i}, i=1,\dots, n-1$,
vanish. For $n=2q+1$ all entries on this diagonal are non-zero for $x_n \neq0$.
Therefore, the principal minor of the matrix $A$, which is obtained by deleting the $n$-th row and $n$-th column, has a non-zero determinant (this is the matrix $\widetilde{A}$).
For $n = 2q+2$ and $x_n \neq 0$ all entries on this ``anti-diagonal'' vanish, except
$\widetilde{A}_{n/2,n/2} = \widetilde{A}_{q+1,q+1}$.
Therefore, the maximal non-degenerate minor is obtained from $A$ by deleting the $n/2$-th and $n$-th rows and columns with the same numbers.
The proposition is proved.

\begin{corollary}
The algebra of Casimir polynomials for $\V_{2q+1}, q=1,2,\dots$, is generated by the linear function
$x_{2q+1}$.
\end{corollary}

For $n=2q+2$ and $x_n \neq 0$ the matrix \eqref{Av} takes the form
$$
A = \begin{pmatrix}
\dots & B & 0 \\
\dots& 0 & 0 \\
C & 0 & 0 \\
0 & 0 & 0 \end{pmatrix},
$$
where $B$ is the $q \times (q+1)$-matrix
$$
B =
\begin{pmatrix}
q x_{q+2} & (q+1) x_{q+3} & \dots & (2q-1) x_{2q+1} & 2q x_{2q+2}  \\
\dots & \dots & \dots & \dots & 0 \\
(q-i) x_{q+1+i} & (q+2-i) x_{q+2+i} & \dots & 0 & 0 \\
\dots & \dots & \dots & 0 & 0 \\
x_{2q+1} & 2 x_{2q+2} & 0 & \dots & 0
\end{pmatrix}
$$
and $C$ is the following $q \times q$-matrix
$$
C =
\begin{pmatrix}
-(q+1)x_{q+3} & -q x_{q+4} & \dots & -3 x_{2q+1} & -2 x_{2q+2} \\
-(q+2) x_{q+4} & -(q+1) x_{q+5} & \dots & -4 x_{2q+2} & 0 \\
\dots & \dots & \dots & 0 & 0 \\
-2q x_{2q+2} & 0 & \dots & 0 & 0
\end{pmatrix}.
$$
The Casimir polynomials must satisfy the equation \eqref{s1}, which can be rewritten as
$$
A v = 0, \ \ v = \begin{pmatrix} \frac{\partial F}{\partial x_1} \\ \frac{\partial F}{\partial x_2} \\ \vdots \\
\frac{\partial F}{\partial x_{2q+2}} \end{pmatrix}.
$$
Obviously, the function $F = x_{2q+2}$ satisfies this system. Let us look for another solution, functionally independent of this one.

Since for $x_{2q+2} \neq 0$ the matrix $C$ is non-degenerate, then
\begin{equation}
\label{null}
\frac{\partial F}{\partial x_1} = \dots = \frac{\partial F}{\partial x_q} = 0.
\end{equation}
Taking into account the above, the system \eqref{s1} reduces to
$$
B w = 0, \ \ w = \begin{pmatrix}\frac{\partial F}{\partial x_{q+1}} \\ \vdots \\ \frac{\partial F}{\partial x_{2q+1}}\end{pmatrix}.
$$
For $x_{2q+2} \neq 0$ the rank of $B$ is equal to $q$. Therefore solutions of the equation $Bw = 0$ form a one-dimensional subspace.

We will extend the matrix $B$ to a square matrix $\widetilde{B}$ by adding a row $(1,0,\dots,0)$:
$$
\widetilde{B} = \begin{pmatrix} & B & \\ 1 & 0 & 0 \end{pmatrix},
$$
and instead of the equation $Bw = 0$ we consider the equation
\begin{equation}
\label{s2}
\widetilde{B} w =
\begin{pmatrix} x_{2q+2}^q \\ 0 \\ \vdots \\ 0\end{pmatrix}.
\end{equation}

This system is solved sequentially as follows:
\begin{equation}
\label{s3}
\frac{\partial F_n}{\partial x_{q+1}} = x^q_{2q+2}, \ \ \ n=2q+2 \geq 4,
\end{equation}
and for $m =0,\dots,q-1$
\begin{equation}
\label{s4}
\frac{\partial F_n}{\partial x_{q+m+2}} = -\frac{1}{(2m+2)x_{2q+2}} \sum_{l=1}^{m+1} (l+m)x_{2q+l-m} \frac{\partial F_n}{\partial x_{q+l}}.
\end{equation}
All quantities $\frac{\partial F_{2q+2}}{\partial x_k}, k=q+1,\dots,2q+1$, are polynomials in
$x_{q+1},\dots$, $x_{2q+2}$.

We will not derive the general form of $F_{2k}$,
but just present the formal expression
\begin{equation}
\label{formal}
F_{2q+2} = \int \left( \frac{\partial F_{2q+2}}{\partial x_{q+1}} dx_{q+1} + \dots +
\frac{\partial F_{2q+2}}{\partial x_{2q+1}} dx_{2q+1}\right. +
\end{equation}
$$
\left. + \frac{\partial F_{2q+2}}{\partial x_{2q+2}} dx_{2q+2} \right),
$$
where $\frac{\partial F_{2q+2}}{\partial x_{q+1}}, \dots$, $\frac{\partial F_{2q+2}}{\partial x_{2q+1}}$
are found from the linear equation \eqref{s2}, and $\Phi = \frac{\partial F_{2q+2}}{\partial x_{2q+2}}$
is determined by the compatibility conditions
$$
\frac{\partial \Phi}{\partial x_i} = \frac{\partial^2 F_{2q+2}}{\partial x_i \, \partial x_{2q+2}},
i=q+1,\dots,2q+1.
$$

Since $\det \widetilde{B} = 2^q \, q! \, x_{2q+2}^q$, all partial derivatives
$\frac{\partial F_{2q+2}}{\partial x_{q+1}}, \dots$, $\frac{\partial F_{2q+2}}{\partial x_{2q+1}}$,
and with them
$\frac{\partial F_{2q+2}}{\partial x_{2q+2}}$, as already noted, are polynomials
in $x_{q+1},\dots,x_{2q+2}$.
This is what the choice of the right-hand side in \eqref{s2} is related to.

We also note the following fact, which immediately follows from \eqref{s3} and \eqref{s4}:

\begin{proposition}
For $n=2q+2$
$$
F_n (x_1,\dots,x_{2q+2}) = \frac{(-1)^q \prod_{m=0}^{q-1} (2m+1)}{2^q\, (q+1)!} x^{q+1}_{2q+1} + x_{2q+2}
G_n (x_{q+1},\dots,x_{2q+2}),
$$
where $G_n$ is a polynomial in $q+2$ variables.
\end{proposition}

We have

\begin{theorem}
The orbits of the coadjoint representation in the coalgebra $\V^\ast_n$

1) for $n=2q+2$ and $x_{2q+2} \neq 0$ are the level surfaces of the polynomials $f_1 = x_{2q+2}$
and $f_2 = F_{2q+2}$, where the polynomial $F_{2q+2}$ is defined by the relations \eqref{s2} and
\eqref{formal} and for $x_{2q+2}=0$ is equal to $\mathrm{const}\cdot x_{2q+1}^{q+1}$;

2) for $n=2q+1$ and $x_{2q+1} \neq 0$ are the hyperplanes
$x_{2q+1} = \mathrm{const} \neq 0$;

3) for $x_n=0$ are the orbits of the coadjoint representation in $\V^\ast_{n-1}$.
\end{theorem}

Note that for $n=2$ the algebra $\V_n$ is commutative and the orbits of the coadjoint representation are points, and for $n=3$ we obtain the Heisenberg algebra and the orbits of the coadjoint representation are hyperplanes $x_3 = \mathrm{const}\neq 0$ or points $(x_1,x_2,0)$.

\begin{corollary}
A typical orbit of the coadjoint representation in the algebra $\V^\ast_n$ has the maximum possible dimension, which is not equal to the dimension of the algebra.
\end{corollary}

Let us present the simplest examples of systems for $F_{2q+2}$ and their solutions:

1) $q=1, n=4$:
$$
\frac{\partial F_4}{\partial x_2}= x_4, \ \frac{\partial F_4}{\partial x_3} = -\frac{1}{2} x_3
$$
and for $F_4$ we may take
$$
F_4(x_2,x_3,x_4) = x_2 x_4 -\frac{1}{4}x_3^2;
$$

2) $q=2, n=6$:
$$
\frac{\partial F_6}{\partial x_3} = x_6^2, \ \
\frac{\partial F_6}{\partial x_4} = -\frac{1}{2} x_5 x_6, \ \
\frac{\partial F_6}{\partial x_5} = \frac{3}{8} x_5^2,
$$
and
$$
F_6(x_3,x_4,x_5,x_6) = x_3 x_6^2 - \frac{1}{2} x_4 x_5 x_6 + \frac{1}{8} x_5^3
$$
up to terms of the form $f(x_6)$;

3) $q=3, n=8$:
$$
\frac{\partial F_8}{\partial x_4}  = x_8^3, \ \ \frac{\partial F_8}{\partial x_5} = -\frac{1}{2} x_7 x_8^2,  \ \
\frac{\partial F_8}{\partial x_6} = -\frac{1}{2} x_6 x_8^2 + \frac{3}{8} x_7^2 x_8,
$$
$$
\frac{\partial F_8}{\partial x_7} = -\frac{1}{2} x_5 x_8^2 + \frac{3}{4} x_6 x_7 x_8 - \frac{15}{48} x_7^3,
$$
and we obtain (up to terms equal to polynomials in $x_8$)
$$
F_8 (x_4,x_5,x_6,x_7,x_8)  = x_4 x_8^3 -\frac{1}{2} x_5x_7 x_8^2 - \frac{1}{4} x_6^2 x_8^2 + \frac{3}{8} x_6 x_7^2 x_8 - \frac{15}{48} x_7^4.
$$

We see that, in contrast to the series of nilpotent algebras considered in the previous section
$\Q_n$, for which the dimensions of typical orbits of the coadjoint representation have the minimum possible dimensions, for the algebras $\V_n$ the situation is the opposite: typical orbits of the coadjoint representation have the maximum possible dimensions.
Therefore, of interest is

\begin{problem}
Do there exist on $\V_n$, with $n \geq 5$, integrable Euler's equations, including those corresponding to geodesic flows of left-invariant metrics on $V_n$?
\end{problem}

\section{The symplectic nilmanifolds $M(2n)$}
\label{secnm}

We recall that the nilpotent Lie algebras $\V_n, n \geq 1$ introduced in \S \ref{secvn}
are defined by the commutation relations satisfied by their generators
$e_1,\dots,e_n = \dim \V_n$:
\begin{equation}
\label{comm1}
[e_i,e_j] = \begin{cases} (j-i)e_{i+j}, & i+j \leq n; \\ 0, & i+j>n. \end{cases}.
\end{equation}

We identify the algebra $\V_n$ with the linear space $\R^n$ with generators $e_1,\dots$, $e_n$ and define the multiplication of vectors from $\R^n$ by the Campbell--Hausdorff formula \eqref{ch}.
This multiplication defines on $\R^n$ the structure of a nilpotent Lie group $V_n$.
The lattice $\Gamma_n = V_n(\Z)$, formed by vectors with integer coordinates,
forms a discrete subgroup in the nilpotent groups $V_n$. The compact nilmanifolds
$$
M(n) = V_n/\Gamma_n,
$$
constructed in \cite{BaT}, have a number of remarkable properties, which were subsequently used and considered in \cite{BaT2,Ba}.

If for each $n$ we denote by $\omega_1,\dots,\omega_n$ the basis of left-invariant $1$-forms dual to the basis $e_1,\dots,e_n$,
then the following relations hold
$$
d\omega_k = (k-2) \omega_1 \wedge \omega_{k-1} + (k-4) \omega_2 \wedge \omega_{k-2} + \dots.
$$
From these it immediately follows that

1) $2$-forms
$$
\Omega_{2n} = (2n - 1)\omega_1 \wedge \omega_{2n} + (2m-3) \omega_2 \wedge_{2m-1} + \dots + \omega_n \wedge \omega_{n+1}
$$
are symplectic forms on the nilmanifolds $M(2n)$;

2) $1$-forms $\omega_{2n+1}$ are contact forms on $M(2n+1)$;

3) there is a natural tower of $S^1$-bundles
$$
\dots \to M(n+1) \to M(n) \to M(n-1) \dots,
$$
where the curvature form of the bundle $M(n+1) \to M(n)$ is
$$
\Omega = (n-1) \omega_1 \wedge \omega_n + (n-3) \omega_2 \wedge \omega_{n-1} + \dots.
$$

In \cite{BaT,BaT2} these nilmanifolds $M(2n)$ were considered as nonformal symplectic manifolds, with the help of which one can construct nonformal simply connected closed symplectic manifolds.
The first examples of such manifolds were constructed in
\cite{BaT98,BaT} for all even dimensions greater than $8$. Since all simply connected closed manifolds of dimension $\leq 6$ are formal, the question remained open in dimension $8$.
Using a method going back to \cite{Guan95}, an example of a simply connected symplectic closed manifold, which is nonformal, was constructed in \cite{FM}.

In \cite{Guan} it was shown that nilpotent algebras
$Q_n$ also generate a tower of nilmanifolds similar to $M(2n)$, and a question posed by its author in the early 1990s was considered: is it true that for a symplectic compact nilmanifold the covering nilpotent group has class  $\leq 2$ as a solvable group. In \cite{Guan} an example is given where the nilpotent group covering a compact symplectic nilmanifold has class $3$ as a solvable group.
The manifolds $M(2n)$ show that such class is generally unbounded. Indeed, the following fact holds.

\begin{theorem}
The nilpotent Lie algebra $\V_n$ has step $k$ as a solvable Lie algebra, where
$$
2^{k}-1 \leq n < 2^{k+1}-1.
$$
\end{theorem}

To prove this theorem, recall that a Lie algebra $\g$
is called solvable if the series of ideals
$$
\g \supset D\g = [\g,\g] \supset D^2\g =
[D^1\g,D^1\g] \supset \dots ... \supset D^i\g = [D^{i-1}\g,D^{i-1}\g] \supset \dots
$$
terminates at the final step. If $D^k\g = 0$ and $D^{k-1}\g \neq 0$, then $k$ is called the step of the solvable algebra $\g$.

It follows from \eqref{comm1} that
$$
\D^1 \V_n = \mathrm{span}\, (e_3,e_4,\dots,e_n), \ \
\D^2 \V_n = \mathrm{span}\, (e_7,e_8,\dots,e_n), \dots,
$$
$$
\D^k \V_n = \mathrm{span}\, (e_{2^{k+1}-1},e_{2^{k+1}},\dots,e_n), \dots .
$$
Theorem 3 follows immediately from these relations.

The classification of compact symplectic nilmanifolds is not yet complete, only answers are known
in low dimensions \cite{BM,Salamon}. This classification is carried out up to finite-sheeted coverings, which in topological language means up to rational homotopy equivalence.

We say that a nilpotent Lie algebra $\g$ is symplectic if on
the corresponding simply connected Lie group $G$ there exists a left-invariant symplectic $2$-form. Maltsev's theorem states that in a simply connected nilpotent Lie group $G$
there exists a discrete subgroup $\Gamma$ such that the quotient space is a closed nilmanifold if and only if with respect to some basis of the Lie algebra
all its structure constants are rational: $c^k_{ij} \in \mathbb{Q}$ \cite{Maltsev}.
Therefore, the classification of closed nilmanifolds (up to rational
homotopy equivalence) reduces to the classification up to isomorphism
of nilpotent Lie algebras over $\mathbb{Q}$. According to \cite{Maltsev}, two nilmanifolds
with isomorphic Lie algebras over $\mathbb{Q}$ are coverings of another nilmanifold.
As I.K.~Babenko recently pointed out to the author, nilmanifolds $M(2n)$ have another important property:

\vskip2mm

{\sl $M(2n), n \geq 2$, are in some sense maximal, since they are not covered by any other nilmanifolds and, in particular, their first homology groups are torsion-free:
$\mathrm{Torsion} \, H_1(M(2n); \Z) = 0$.}

\vskip2mm

It is known that simply connected homogeneous symplectic manifolds are symplectomorphic to orbits of coadjoint representations of compact semisimple Lie groups. In \cite{Efimov} it was proved that the magnetic geodesic flow on any such orbit, where the metric is given by the Killing form and the magnetic field is given by the Kirillov form (i.e. generated by the Lie--Poisson brackets), is (non-commutatively) integrable. It was thus shown that on any simply connected homogeneous symplectic manifold there exists a metric such that in the magnetic field given by the symplectic form the magnetic geodesic flow is integrable. It is natural to ask whether there is an analogue of this result for non-simply connected homogeneous symplectic manifolds. In particular, we have

\begin{problem}
Do there exist Riemannian metrics (including left-invariant ones) on nilmanifolds $M(2n)$ such that the magnetic geodesic flows corresponding to them and the magnetic fields $\Omega_{2n}$ are integrable?
\end{problem}

In \cite{Guan} it is stated that nilmanifolds corresponding to groups $Q_{2n}$ are symplectic for all $n$. The corresponding Lie algebras form a tower of central extensions
$$
\Q_1 \subset \Q_2 \subset \dots \subset \Q_n \subset \dots,
$$
analogously to the algebras
$$
\V_1 \subset \V_2 \subset \dots \subset \V_n \subset \dots.
$$
The first four members of these towers are isomorphic:
$$
\Q_k = \V_k \ \ \mbox{for} \ \ k=1, 2, 3, 4,
$$
and already $\Q_k \neq \V_k$ for $k=5,6$.
It is natural to

\begin{problem}
Is it true that any sequence of nontrivial central extensions
of Lie algebras over $\mathbb{Q}$
$$
\mathcal{A}_1 \subset \mathcal{A}_2 \subset \dots \subset \mathcal{A}_k, \ \ \dim \mathcal{A}_i = i,
$$
such that its even-dimensional members $\mathcal{A}_{2n}$ are symplectic, can be extended indefinitely to a sequence with the same properties?
If such extensions exist for some finite sequence, how can they be classified?
\end{problem}

Since odd-dimensional nilmanifolds $V_{2n+1}$ are contact, this question can be posed for contact Lie algebras as well.

The Lie algebras $\V_n$ and nilmanifolds $M(2n)$ have other interesting properties:

1) the Lie algebras $\Q_n$ and $\V_n$ are filiform, i.e., they are obtained from a one-dimensional algebra by successive one-dimensional central extensions \cite{Vergne}.
It follows from the commutation relations that they are also $\N$-graded. As proved in \cite{Mill}, any symplectic filiform Lie algebra of dimension $n \geq 12$ is a
$\N$-graded deformation of either $\Q_n$ or $\V_n$. In \cite{Mill} the moduli space of such deformations of the algebras $\V_n, n \geq 16$ was also described;

2) soon after the appearance of the preprint of the work \cite{BaT}, Buchstaber pointed out an interesting connection between the nilmanifolds $M(2n)$ constructed in it and the Land\-we\-be\-ra--Novikov algebras \cite{Buch}.

\end{document}